\documentclass[reqno,a4paper,11pt]{amsart}

\usepackage{amsmath,amstext,amsfonts,amsbsy,eucal,amssymb}
\usepackage[latin1]{inputenc}

\numberwithin{equation}{section}

\newtheorem{theorem}{Theorem}[section]

\newtheorem{proposition}[theorem]{Proposition}

\newtheorem{conjecture}[theorem]{Conjecture}

\theoremstyle{definition}

\newtheorem{remark}[theorem]{Remark}

\begin{document}

\parskip 4pt
\baselineskip 16pt


\title[Conjectures involving a generalization of the sums of powers of integers]
{Conjectures involving a generalization of the sums of powers of integers}

\author[Andrei K. Svinin]{Andrei K. Svinin}
\address{Andrei K. Svinin, 
Matrosov Institute for System Dynamics and Control Theory of 
Siberian Branch of Russian Academy of Sciences,
P.O. Box 292, 664033 Irkutsk, Russia}
\email{svinin@icc.ru}
%
%

\date{\today}



\begin{abstract}
Some class of  sums which naturally include the sums of powers of integers is considered. A number of conjectures concerning a representation of these sums is made.
\end{abstract}

\maketitle

\section{Introduction}

It is common  knowledge  that the sum of powers of integers
\begin{equation}
S_m(n):=\sum_{q=1}^nq^m
\label{014}
\end{equation}
is a polynomial $\sigma_m(n)$ in $n$ of degree $m+1$. This proposition can be showed  by using Pascal's elementary proof (see, for example \cite{Beardon}).  Polynomials $\sigma_m(n)$ for any integer $m\geq 1$ are defined  in terms of the Bernoulli numbers: 
\[
\sigma_m(n)=\frac{1}{m+1}\sum_{q=0}^m(-1)^q{m+1\choose q}B_{q}n^{m-q+1}.
\] 
The numbers $B_k$, in turn, are defined by a recursion
\[
\sum_{q=0}^{k}{k+1\choose q}B_q=0,\;\;
B_0=1.
\]
They were discovered by Bernoulli in 1713 but as explained in \cite{Edwards1},  they were known by Faulhaber before this time. 

Faulhaber's theorem says that for any odd $m\geq 3$, the sum (\ref{014}) is expressed as a polynomial of $S_1(n)$. It is common knowledge that $S_1(n)=n(n+1)/2$. Let us denote $t:=n(n+1)$. Therefore, the Faulhaber's theorem asserts that for any $k\geq 1$, $S_{2k+1}(n)$ can be expressed as a polynomial in $\mathbb{Q}[t]$. These polynomials are referred as Faulhaber's polynomials. Also, it is worth remembering that it was shown by Jacobi  in \cite{Jacobi} that $S_{2k}(n)$ is expressed as a polynomial in $(2n+1)\mathbb{Q}[t]$. 

There exists a variety of different modifications and generalizations of  Faulhaber's theorem in the literature. Faulhaber himself considered $r$-fold sums $S_{m}^{r}(n)$ which are successively defined by
\[
S_{m}^{r}(n)=\sum_{q=1}^{n}S_{m}^{r-1}(q),\;\;
r\geq 1
\]  
beginning from $S_{m}^{0}(n):=n^m$. He observed that $S_{m}^{r}(n)$ can be expressed as a polynomial in $n(n+r)$ if $m-r$ is even \cite{Knuth}.

In this paper we consider some class of sums $S_{k, j}(n)$ which include classical sums of powers of integers (\ref{014}) with odd $m=2k+1$.  These sums will be introduced in section \ref{s2} and \ref{s3}.  In the particular case $n=1$, they represent some class of binomial sums studied in \cite{Tuenter}. Our crucial idea comes from this work. We  briefly describe the results of this article in section \ref{s4}. In section \ref{s5}, we formulate our main conjecture concerning a representation of the sums  $S_{k, j}(n)$ in terms of some polynomials whose coefficients rationally depend on $n$. In section \ref{s6}, we discuss the relationship of this representation with a possible generalization of the Faulhaber's theorem. In section \ref{s7}, we show the explicit form of some $n$-dependent coefficients of these polynomials and speculate about the relationship of these coefficients with the number of $q$-points in the simple symmetric $(2k+1)$-Venn diagram. 

\section{The sums $s_{k, j}(n)$}
\label{s2}

Let us define the numbers $\{C_{j, r}(n) : j\geq 1, r=0,\ldots, j(n-1)\}$ through the relation
\begin{equation}
\left(\sum_{q=0}^{n-1}a^{q+1}\right)^j=\sum_{q=0}^{j(n-1)}C_{j, q}(n)a^{q+j},
\label{013}
\end{equation}
where $a$ is supposed to be an auxiliary variable, taking into account that $a^la^m=a^{l+m}$. For example, in the case $n=1$, we have only $C_{j, 0}(1)=1$, while  in the case $n=2$, we  get
\[
\left(a+a^2\right)^j=\sum_{q=0}^j{j\choose q}a^{j+q},
\]
and therefore  $C_{j, q}(2)$  is a binomial coefficient.  From
\begin{eqnarray}
\left(\sum_{q=0}^{n-1}a^{q+1}\right)^j&=&\left(\sum_{q=0}^{(j-1)(n-1)}C_{j-1, q}(n)a^{q+j-1}\right)\left(\sum_{q=0}^{n-1}a^{q+1}\right)\nonumber\\
&=&C_{j-1, 0}(n)a^{j}+\left(C_{j-1, 0}(n)+C_{j-1, 1}(n)\right)a^{j+1}+\cdots+\left(C_{j-1, 0}(n)+\cdots\right.\nonumber\\
&&\left.+C_{j-1, n-1}(n)\right)a^{j+n-1}\nonumber\\
&&+\left(C_{j-1, 1}(n)+\cdots+C_{j-1, n)}(n)\right)a^{j+n}+\cdots+C_{j-1, (j-1)(n-1)}(n)a^{jn}\nonumber
\end{eqnarray}
we deduce  that
\begin{equation}
C_{j,r}(n)=\sum_{q=r-n+1}^rC_{j-1, q}(n)
\label{008}
\end{equation}
in assumption that $C_{j-1, q}(n)=0$ for $q<0$ and $q>(j-1)(n-1)$. In the case $n=2$, the relation (\ref{008}) becomes a well-known property for the binomial coefficients:
\[
{j\choose r}={j-1\choose r-1}+{j-1\choose r}.
\]

Clearly, the coefficient $C_{j,r}(n)$ can be represented as a proper sum of multinomial coefficients.
Putting $a=1$ into (\ref{013}), we get the property
\[
\sum_{q=0}^{j(n-1)}C_{j, q}(n)=n^j
\]
for these numbers which evidently generalize the well-known attribute of binomial coefficients. Also it is easy to get the following property:
\[
\sum_{q=0}^{j(n-1)}qC_{j, q}(n)=\frac{j(n-1)}{2}n^j.
\]

With the coefficients $C_{j, q}(n)$ we define the sum $s_{k, j}(n)$ by
\[
s_{k, j}(n):=\sum_{q=0}^{j(n-1)}C_{j, q}(n)x_{j+q},
\]
where $x_r:=r^{2k+1}$. Let us remark that $s_{k, 1}(n)=S_{2k+1}(n)$. Next, we would like to give a suitable expression for the sums $s_{k, 2}(n)$, which will be used in the following section.
It is evident that
\[
C_{2, q}(n)=\left\{
\begin{array}{l}
q+1,\;\; q=0,\ldots, n-2, \\[0.3cm]
2n-q-1,\;\; q=n-1,\ldots, 2n-2 
\end{array}
\right.
\]
and therefore
\begin{eqnarray}
s_{k, 2}(n)&:=&\sum_{q=0}^{2n-2}C_{2, q}(n)x_{q+2}\nonumber\\
&=&\sum_{q=0}^{n-2}(q+1)x_{q+2}+\sum_{q=n-1}^{2n-2}(2n-q-1)x_{q+2}.\nonumber
\end{eqnarray}
Shifting $q\rightarrow q-2$, we get
\[
s_{k, 2}(n)=\sum_{q=1}^{n}(q-1)x_{q}+\sum_{q=n+1}^{2n}(2n-q+1)x_{q},
\]
and since $q-1=(2n-q+1)-(2n-2q+2)$, we deduce that
\begin{equation}
s_{k, 2}(n)=-\sum_{q=1}^{n}(2n-2q+2)x_{q}+\sum_{q=1}^{2n}(2n-q+1)x_{q}.
\label{020}
\end{equation}

\section{The sums $S_{k, j}(n)$ and $\tilde{S}_{k, j}(n)$}
\label{s3}

Let us define a sum
\[
\tilde{S}_{k, j}(n):=\sum_{\{\lambda\}\in B_{j, jn}}\left\{\lambda_1^{2k+1}+(\lambda_2-n)^{2k+1}+\cdots+(\lambda_j-jn+n)^{2k+1}\right\},
\]
where $B_{j, jn}:=\{\lambda_k : 1\leq \lambda_1\leq \cdots \leq \lambda_j\leq jn\}$ and
\begin{equation}
S_{k, j}(n)=\sum_{q=0}^{j-1}{j(n+1)\choose q} s_{k, j-q}(n).
\label{003}
\end{equation}
It is obvious that 
\[
\tilde{S}_{k, 1}(n)=S_{k, 1}(n)=s_{k, 1}(n)=S_{2k+1}(n).
\]

\begin{conjecture}
\label{con1}
For any $k\geq 0$ and $j\geq 1$, one has
\begin{equation}
\tilde{S}_{k, j}(n)=S_{k, j}(n).
\label{001}
\end{equation}
\end{conjecture}
Let us prove (\ref{001}) in the particular case $j=2$. Using simple arguments we deduce that 
\begin{eqnarray}
\tilde{S}_{k, 2}(n)&=&\sum_{\{\lambda\}\in B_{2, 2n}}\left\{\lambda_1^{2k+1}+(\lambda_2-n)^{2k+1}\right\}\nonumber\\
&=&\sum_{q=1}^{2n}(2n-q+1)x_q+\sum_{q=1}^{2n}(2n-q+1)x_{n-q+1}.
\label{002}
\end{eqnarray}
Notice that  we have several $x_r$'s with negative subscript $r$ in (\ref{002}). Clearly, we must put $x_r=-x_{-r}$. Taking into account this rule, we rewrite the second sum in (\ref{002})  as
\[
\sum_{q=1}^{n}(n+q)x_{q}-\sum_{q=1}^{n}(n-q)x_{q}=2\sum_{q=1}^{n}qx_{q}.
\]
Adding to (\ref{002}) the terms
\[
-\sum_{q=1}^{n}(2n-2q+2)x_{q}+\sum_{q=1}^{n}(2n-2q+2)x_{q}
\]
and taking into account (\ref{020}), we finally get
\begin{eqnarray}
\tilde{S}_{k, 2}(n)&=&s_{k, 2}(n)+(2n+2)\sum_{q=1}^{n}x_q\nonumber\\
&=&s_{k, 2}(n)+{2(n+1)\choose 1}s_{k, 1}(n)\nonumber\\
&=&S_{k, 2}(n).\nonumber
\end{eqnarray}
Therefore, we proved conjecture \ref{con1} in the particular case $j=2$.

\section{Special case of the sums $S_{k, j}(n)$}
\label{s4}

In the case $n=1$, the sum (\ref{003}) becomes
\begin{equation}
S_{k, j}(1)=\sum_{q=0}^{j-1}{2j\choose q} (j-q)^{2k+1}.
\label{004}
\end{equation}
This type of the sums was studied in \cite{Bruckman}, \cite{Strazdins}, \cite{Tuenter}. To be precise, the authors of these works considered binomial sums of the form
\[
\mathcal{S}_{m}(j)=\sum_{q=0}^{2j}
{2j\choose q}|j-q|^m.
\]
It is evident that $S_{k, j}(1)=\mathcal{S}_{2k+1}(j)/2$.

Some bibliographical remarks are as follows.  Bruckman  in \cite{Bruckman} asked to prove that $\mathcal{S}_3(j)=j^2\left(
\begin{array}{c}
2j\\
j
\end{array}
\right)$. Strazdins in \cite{Strazdins} solved this problem and conjectured that $\mathcal{S}_{2k+1}(j)=\tilde{P}_{k}(x)|_{x=j}{2j\choose j}$ with some monic polynomial $\tilde{P}_{k}(x)$ for any $k\geq 0$. Tuenter showed in \cite{Tuenter} that it is almost true. More exactly, he proved that there exists a sequence of polynomials $P_k(x)$ such that
\[
\mathcal{S}_{2k+1}(j)=P_{k}(x)|_{x=j}j
{2j\choose j}=P_{k}(x)|_{x=j}\frac{(2j)!}{(j-1)!j!}.
\]
One can see that polynomial $\tilde{P}_{k}(x)$ is monic only for $k=0, 1$. The following proposition \cite{Tuenter} is verified by direct computations.
\begin{proposition}
The sums $S_{k, j}(1)$ satisfy the recurrence relation
\begin{equation}
S_{k, j}(1)=j^2S_{k-1, j}(1)-2j(2j-1)S_{k-1, j-1}(1).
\label{rec1}
\end{equation}
\end{proposition}
\begin{remark}
This was, in fact, proved  for the sums $\mathcal{S}_{2k+1}(j)$ in \cite{Tuenter}.
\end{remark}

Let us introduce the sequence of positive integer numbers:
\[
g_1=1,\;\;g_j=\frac{j}{(j-1)!}\prod_{q=1}^{j-1}(j+q),
\;\;j\geq 2.
\]
As can be checked, this sequence  satisfies the recurrence relation
\begin{equation}
g_{j+1}=2\frac{2j+1}{j}g_j.
\label{016}
\end{equation}
Some of the first numbers $g_j$ are given in the following table:
\begin{center}
\begin{tabular}{|c|c|c|c|c|c|c|c|c|c|}
\hline
$j$&$1$&$2$&$3$&$4$&$5$&$6$&$7$&$8$&9  \\
\hline
$g_{j}$&$1$&6&30&140& 630&2772 &12012 &51480 &218790\\
\hline
\end{tabular}
\end{center}
They are included in the \textrm{A002457} integer sequence in \cite{Sloane}.

Let us write 
\[
S_{k, j}(1)=P_{k, j}g_j
\]
with some numbers  $P_{k, j}$ to be calculated. Taking into account (\ref{016}), it can be easily seen that (\ref{rec1}) is valid if the relation
\begin{equation}
P_{k+1, j}=j^2\left(P_{k,j}-P_{k,j-1}\right)+jP_{k, j-1}
\label{007}
\end{equation}
holds. Now we need the polynomials studied in \cite{Tuenter} which are defined by a recurrence relation 
\begin{equation}
P_{k+1}(x)=x^2\left(P_k(x)-P_{k}(x-1)\right)+xP_k(x-1)
\label{006}
\end{equation}
with initial condition $P_0(x)=1$. The first eight polynomials yielded by (\ref{006}) are as follows:
\[
P_0(x)=1,\;\;
P_1(x)=x,\;\;
P_2(x)=x(2x-1),
\]
\[
P_3(x)=x(6x^2-8x+3),
\]
\[
P_4(x)=x(24x^3-60x^2+54x-17),
\]
\[
P_5(x)=x(120x^4-480x^3+762x^2-556x+155),
\]
\[
P_6(x)= x(720x^5-4200x^4+10248x^3-12840x^2+8146x-2073).
\]
\begin{eqnarray}
P_7(x)&=& x\left(5040x^6-40320x^5+139440x^4-263040x^3\right. \nonumber\\
&&\left.+282078x^2-161424x+38227\right).\nonumber
\end{eqnarray}

Introducing a recursion operator $R:=x^2\left(1-\Lambda^{-1}\right)+x\Lambda^{-1}$, where $\Lambda$ is a shift operator acting as $\Lambda(f(x))=f(x+1)$, one can write $P_k(x)=R^k(1)$. One could notice that the polynomial $P_k(x)$ has $k!$ as a coefficient at $x^k$. In addition, as was noticed in \cite{Tuenter}, the constant terms of the polynomials $P_k(x)/x$ form a sequence of the Genocchi numbers with opposite sign:
\[
-G_2=1,\;\;
-G_4=-1,\;\;
-G_6=3,\;\;
-G_8=-17,\;\;
-G_{10}=155,\;\;
-G_{12}=-2073,\ldots
\]
Recall that the Genocchi numbers are defined in terms of the generating function
\[
\frac{2x}{e^x+1}=\sum_{q\geq 1}G_{q}\frac{x^q}{q!}
\]
and are related to the Bernoulli numbers as $G_{2k}=2\left(1-2^{2k}\right)B_{2k}$.

Comparing (\ref{006}) with (\ref{007}) we conclude that $P_{k, j}=P_k(x)|_{x=j}$. 
\begin{remark}
As was noticed in \cite{Tuenter}, the polynomials $P_{k}(x)$ can be obtained as a special case of Dumont-Foata polynomials of three variables \cite{Dumont}.
\end{remark}
\section{General case of the sums $S_{k, j}(n)$}
\label{s5}

Let us remember that $t:=n(n+1)$. Our main conjecture is as follows.
\begin{conjecture}
There exists some polynomials $P_k(t, x)$ in $x$ of degree $k$ whose coefficients rationally depend on $t$ such that 
\[
S_{k,j}(n)=\frac{t^{k+1}}{2^{k+1}}P_k(t, x)|_{x=j}g_j(n),
\]
where
\[
g_1(n)=1,\;\;
g_j(n):=\frac{j}{(j-1)!}\prod_{q=1}^{j-1}(jn+q).
\]
\end{conjecture}
The first eight polynomials $P_k(t, x)$ are 
\[
P_0(t, x)=1,\;\;
P_1(t, x)=x,\;\;
P_2(t, x)=x\left(2x-\frac{2(t+1)}{3t}\right),\;\;
\]
\[
P_3(t, x)=x\left(6x^2-\frac{16(t+1)}{3t}x+\frac{4(t+1)^2}{3t^2}\right),\;\; 
\]
\[
P_4(t, x)=x\left(24x^3-\frac{40(t+1)}{t}x^2+\frac{24(t+1)^2}{t^2}x-\frac{24(t+1)^3+8t^2}{5t^3}\right),
\]
\begin{eqnarray}
P_5(t, x)&=& x\left(120x^4-\frac{320(t+1)}{t}x^3+\frac{1016(t+1)^2}{3t^2}x^2\right.\nonumber\\
&&\left.-\frac{160(t+1)^3+32t^2}{t^3}x+\frac{80\left((t+1)^4+t^2(t+1)\right)}{3t^4}\right),\nonumber
\end{eqnarray}
\begin{eqnarray}
P_6(t, x)&=& x\left(720x^5-\frac{2800(t+1)}{t}x^4+\frac{13664(t+1)^2}{3t^2}x^3\right.\nonumber\\
&&-\frac{55936(t+1)^3+7632t^2}{15t^3}x^2+\frac{22112(t+1)^4+13664t^2(t+1)}{15t^4}x\nonumber\\
&&\left.-\frac{22112(t+1)^5+44224t^2(t+1)^2}{105t^5}\right),\nonumber
\end{eqnarray}
\begin{eqnarray}
P_7(t, x)&=&x\left(5040x^6-\frac{26880(t+1)}{t}x^5+\frac{185920(t+1)^2}{3t^2}x^4\right.\nonumber\\
&&-\frac{76800(t+1)^3+7680t^2}{t^3}x^3\nonumber\\
&&+\frac{157088(t+1)^4+67968t^2(t+1)}{3t^4}x^2\nonumber\\
&&-\frac{17920(t+1)^5+22528t^2(t+1)^2}{t^5}x\nonumber\\
&&\left. +\frac{6720(t+1)^6+22400t^2(t+1)^3+1344t^4}{3t^6}\right).\nonumber
\end{eqnarray}
Unfortunately, in the general case we do not know a recursion relation for these polynomials except for the case $t=2$. One can check that $P_k(t, x)|_{t=2}=P_k(x)$, where $P_k(x)$ are Tuenter's polynomials introduced above.

It can be noticed that the  coefficients of the  polynomials $P_k(t, x)$ have a special form. Namely, let 
\begin{equation}
P_k(t, x)=x\left(p_{k, 0}(t)x^{k-1}-p_{k, 1}(t)x^{k-2}+\cdots+(-1)^{k-1}p_{k, k-1}(t)\right).
\label{005}
\end{equation} 
Based on actual calculations, we are led to the following.
\begin{conjecture} \label{4.2} 
The coefficients of $t$-dependent polynomial (\ref{005}) are given by 
\[
p_{k, j}(t)=\frac{r_{k, j}(t)}{t^j},
\]
where the polynomials $r_{k, j}(t)$ are of the form
\[
r_{k, j}(t)=\sum_{q=0}^{m}\alpha_{k, j, q}t^{2q}(t+1)^{j-3q}
\]
with rational positive  nonzero numbers $\alpha_{k, j, q}$. Here $m\geq 0$, by definition,  is the result of division of the number $j$ 
by $3$ with some remainder $l$, that is, $j=3m+l$. 
\end{conjecture}
It should be noticed that the last conjecture is quite strong. With this conjecture the number $N_k$ of parameters which entirely define the polynomial $P_{k}(t, x)$ for small values of $k$  is given in the following table:
\begin{center}
\begin{tabular}{|c|c|c|c|c|c|c|c|c|c|}
\hline
$k$&$1$&$2$&$3$&$4$&$5$&$6$&$7$&$8$&9  \\
\hline
$N_{k}$&$1$&2&3&5&7&9&12&15 &18\\
\hline
\end{tabular}
\end{center}
These numbers form the \textrm{A001840} integer sequence in \cite{Sloane}. We have, in fact, verified  conjecture \ref{4.2}  up to $k=11$ and calculated all corresponding coefficients $\alpha_{k, j, q}$.

It is evident that for all the polynomials $P_k(t, x)$ we have $p_{k, 0}(t)=k!$. Actual calculations show that
\[
\alpha_{k, 1, 0}=\frac{(k-1)(k+1)}{9}k!,\;\;k\geq 2,
\]
\[
\alpha_{k, 2, 0}=\frac{(k-2)(k+1)(5k^2+k-3)}{810}k!,\;\;k\geq 3,
\]
\[
\alpha_{k, 3, 0}= \frac{(k-3)(k+1)(175k^4-70k^3-724k^2+643k-690)}{765450}k!,\;\;k\geq 4,
\]
\[
\alpha_{k, 3, 1}= \frac{(k-3)(k+1)(2k^2-4k+5)}{1575}k!,\;\;k\geq 4
\]
etc. Looking at these patterns we may guess that
\[
\alpha_{k, j, 0}=(k-j)(k+1)p_j(k)k!,\;\;k\geq j+1
\]
with some polynomial $p_j(k)$ of degree $2j-2$.

\section{Faulhaber's theorem}
\label{s6}

With  the polynomials  $P_k(t, x)$  we are able, for example, to calculate
\[
S_{0,j}(n)=\frac{t}{2}\frac{j}{(j-1)!}\prod_{q=1}^{j-1}(jn+q),\;\;
S_{1,j}(n)=\frac{t^2}{4}\frac{j^2}{(j-1)!}\prod_{q=1}^{j-1}(jn+q),
\]
\[
S_{2,j}(n)=\frac{t^2}{12}\left\{(3j-1)t-1\right\}\frac{j^2}{(j-1)!}\prod_{q=1}^{j-1}(jn+q),
\]
\[
S_{3,j}(n)=\frac{t^2}{24}\left\{(9j^2-8j+2)t^2-(8j-4)t+2\right\}\frac{j^2}{(j-1)!}\prod_{q=1}^{j-1}(jn+q),
\]
\begin{eqnarray}
S_{4,j}(n)&=&\frac{t^2}{20}\left\{(15j^3-25j^2+15j-3)t^3-(25j^2-30j+10)t^2\right.\nonumber\\
           &&\left.+(15j-9)t-3\right\}\frac{j^2}{(j-1)!}\prod_{q=1}^{j-1}(jn+q),\nonumber
\end{eqnarray}
\begin{eqnarray}
S_{5,j}(n)&=&\frac{t^2}{24}\left\{(45j^4-120j^3+127j^2-60j+10)t^4\right.\nonumber\\
&&-(120j^3-254j^2+192j-50)t^3+(127j^2-180j+70)t^2\nonumber\\
&&\left.-(60j-40)t+10\right\}\frac{j^2}{(j-1)!}\prod_{q=1}^{j-1}(jn+q),\nonumber
\end{eqnarray}
\begin{eqnarray}
S_{6,j}(n)&=&\frac{t^2}{840}\left\{(4725j^5-18375j^4+29890j^3-24472j^2+9674j-1382)t^5\right.\nonumber\\
&&-(18375j^4-59780j^3+76755j^2-44674j+9674)t^4\nonumber\\
&&+(29890j^3-73416j^2+64022j-19348)t^3-(24472j^2+38696j-16584)t^2\nonumber\\
&&\left.+(9674j-6910)t-1382\right\}\frac{j^2}{(j-1)!}\prod_{q=1}^{j-1}(jn+q).\nonumber
\end{eqnarray}
\begin{eqnarray}
S_{7,j}(n)&=&\frac{t^2}{48}\left\{(945j^6-5040j^5+11620j^4-14400j^3+9818j^2-3360j+420)t^6\right.\nonumber\\
&&-(5040j^5-23240j^4+44640j^3-43520j^2+21024j-3920)t^5\nonumber\\
&&+(11620j^4-43200j^3+63156j^2-42048j+10584)t^4\nonumber\\
&&-(14400j^3-39272j^2+37824j-12600)t^3\nonumber\\
&&+(9818j^2-16800j+7700)t^2\nonumber\\
&&\left.-(3360j-2520)t+420\right\}\frac{j^2}{(j-1)!}\prod_{q=1}^{j-1}(jn+q).\nonumber
\end{eqnarray}
Looking at these patterns we  suggest the following.
\begin{conjecture}
\label{4.3}
$S_{k,j}(n)$ is expressed as a polynomial in $\prod_{q=1}^{j-1}(jn+q)\mathbb{Q}[t]$.
\end{conjecture}

For $j=1$,  conjecture \ref{4.3} becomes the well-known  Faulhaber's theorem \cite{Faulhaber} which was, in fact, proved by Jacobi in  \cite{Jacobi}. The first eight  Faulhaber's polynomials  are as follows:
\[
S_{0, 1}(n)=\frac{t}{2},\;\;
S_{1, 1}(n)=\frac{t^2}{4},\;\;
S_{2, 1}(n)=\frac{t^2}{12}(2t-1),\;\;
S_{3, 1}(n)=\frac{t^2}{24}\left(3t^2-4t+2\right),
\]
\[
S_{4, 1}(n)=\frac{t^2}{20}\left(2t^3-5t^2+6t-3\right),\;\;
S_{5, 1}(n)=\frac{t^2}{24}\left(2t^4-8t^3+17t^2-20t+10\right),
\]
\[
S_{6, 1}(n)= \frac{t^2}{840}\left(60t^5-350t^4+1148t^3-46584t^2+2764t-1382\right),
\]
\[
S_{7, 1}(n)= \frac{t^2}{48}t^2(3t^6-24t^5+112t^4-352t^3+718t^2-840t+420).
\]
One usually writes
\[
S_{k, 1}(n)=\frac{1}{2(k+1)}\sum_{q=0}^{k}A_{q}^{(k+1)}t^{k-q+1},
\]
where $A^{(k)}_0=1$ and $A^{(k)}_{k-1}=0$.  We know quite a lot about the coefficients $A_q^{(k)}$. Jacobi proved that the coefficients $A_{q}^{(k)}$ satisfy the recurrence relation
\[
(2k+2)(2k+1)A_q^{(k)}=2(k-q+1)(2k-2q+1)A_{q}^{(k+1)}+(k-q+1)(k-q+2)A_{q-1}^{(k+1)}
\]
and tabulated some of them. It was shown by Knuth in \cite{Knuth} that these coefficients satisfy the quite simple implicit recurrence relation
\begin{equation}
\sum_{q=0}^r{k-q\choose 2r+1-2q}A_{q}^{(k)}=0,\;\;
r>0
\label{009}
\end{equation}
which yields an infinite triangle system of equations from which one easily obtains
\[
A_{1}^{(k)}=-\frac{(k-2)k}{6},\;\;
A_{2}^{(k)}=\frac{(k-3)(k-1)k(7k-8)}{360},
\]
\[
A_{3}^{(k)}=-\frac{(k-4)(k-2)(k-1)k(31k^2-89k+48)}{15120},
\]
\[
A_{4}^{(k)}= \frac{(k-5)(k-3)(k-2)(k-1)k(127k^3-691k^2+1038k-384)}{604800}
\]
etc. Gessel and Viennot  showed in \cite{Gessel} that the solution of system (\ref{009}) can be represented as a $k\times k$ determinant
\[
A_{q}^{(k)}=\frac{1}{(1-k)\cdots(q-k)}\left|
\begin{array}{ccccc}
{k-q+1 \choose 3}    &{k-q+1 \choose 1}      &0     & \cdots   &    0  \\[0.1cm]
{k-q+2 \choose 5}    &{k-q+2 \choose 3}      &{k-q+2 \choose 1}     & \cdots   &    0  \\[0.1cm]
\vdots               &\vdots                 &\vdots          &          &    \vdots  \\[0.1cm]
{k-1 \choose 2k-1}    &{k-1 \choose 2k-3}      &{k-1 \choose 2k-5}     & \cdots   &  {k-1 \choose 1} \\[0.1cm]
{k \choose 2k+1}    &{k-1 \choose 2k-1}      &{k-1 \choose 2k-3}     & \cdots   &  {k \choose 3} 
\end{array}
\right|
\]
(see also \cite{Edwards}). A bivariate generating function for the coefficients $A_{q}^{(k)}$ was obtained in \cite{Gessel1}.

\section{The conjectural relationship of the coefficients $p_{k, k-1}(t)$ to simple symmetric Venn diagrams}
\label{s7}

Let us rewrite the coefficients $p_{k, k-1}(t)$  being expressed via $n$, that is,
\[
p_{k, k-1}(n)=p_{k, k-1}(t)|_{t=n(n+1)}.
\]
For example, 
\[
p_{1, 0}(n)=1,\;\;
p_{2, 1}(n)= \frac{2}{3}\frac{n^2+n+1}{n(n+1)},\;\;
p_{3, 2}(n)= \frac{4}{3}\frac{n^4+2n^3+3n^2+2n+1}{n^2(n+1)^2},
\]
\[
p_{4, 3}(n) = \frac{24}{5}\frac{n^6+3n^5+\frac{19}{3}n^4+\frac{23}{3}n^3+\frac{19}{3}n^2+3n+1}{n^3(n+1)^3},
\]
\[
p_{5, 4}(n)= \frac{80}{3}\frac{n^8+4n^7+11n^6+19n^5+23n^4+19n^3+11n^2+4n+1}{n^4(n+1)^4},
\]
\[
p_{6, 5}(n)= \frac{22112}{105}\frac{n^{10}+5n^{9}+17n^{8}+38n^7+61n^6+71n^5+61n^4+38n^3+17n^2+5n+1}{n^5(n+1)^5}
\]
etc. Looking at these patterns we see that 
\begin{equation}
p_{k, k-1}(n)=c_k\frac{v_k(n)}{n^{k-1}(n+1)^{k-1}},
\label{012}
\end{equation}
where $v_k(n)$ is a monic polynomial of degree $2k-2$. All these polynomials are invariant with respect to the transformation
\begin{equation}
v_k(n)\mapsto n^{2k-2}v_k\left(\frac{1}{n}\right).
\label{015}
\end{equation}
Also it is worth to remark that the polynomial $v_4(n)$ unlike the others, has several fractional coefficients.
\begin{conjecture}
The polynomials $v_k(n)$ are given by
\begin{equation}
v_k(n)=\sum_{q=1}^{2k-1}\frac{{2k\choose q}+(-1)^{q+1}}{2k+1}n^{2k-q-1},
\label{010}
\end{equation}
while the coefficients $c_k$ are expressed via the Bernoulli numbers as
\begin{equation}
c_{k}=(-1)^{k+1}(2k+1)2^kB_{2k}.
\label{011}
\end{equation}
\end{conjecture}
Let us notice that if (\ref{010}) is valid then the invariance of the corresponding polynomial  with respect to (\ref{011}) is obvious in virtue of the invariance of binomial coefficients.  

Let $p=2k+1$. It is known that if $p$ is prime then  
\[
T(p,q)=\frac{{p-1\choose q}+(-1)^{q+1}}{p},\;\;
p\geq 5
\]
is the number of $q$-points on the left side of a crosscut of the simple symmetric $p$-Venn diagram  \cite{Mamakani}. This integer sequence is known as the \textrm{A219539} sequence in \cite{Sloane}. 
It is evident that the  row sum 
\[
t_p:=\sum_{q=1}^{p-2}T(p,q)=\frac{2^{p-1}-1}{p}
\]
can be calculated as $v_k(n)|_{n=1}$. The Fermat quotients $(2^{p-1}-1)/p$ for prime $p$'s form the integer sequence \textrm{A007663}  in \cite{Sloane}. Taking into account (\ref{012}) and (\ref{011}), we get
\[
p_{k, k-1}(n)|_{n=1}=2(2^{2k}-1)B_{2k}=-G_{2k}.
\]

\section{Discussion}

In the paper we have considered some class of sums $S_{k, j}(n)$ and conjectured a representation of these sums  in terms of a sequence of the polynomials $\{P_k(t, x) : k\geq 0\}$. This assumption  resulted from computational experiments and is supported  by a large amount of actual calculations. For $n=1$, we get the well-known results from \cite{Tuenter}. This also confirms our assumptions. The conjectural relationship of several coefficients of the polynomials $P_k(t, x)$ being expressed via $n$ to simple symmetric Venn diagrams is quite unexpected and requires explanation.

\section*{Acknowledgments}
I wish to thank the referee for carefully reading the manuscript and for remarks which enabled the presentation of the paper to be improved.
This  work  was  supported  in  part  by the Council for Grants of the
 President  of  Russian  Foundation  for  state  support of the leading
 scientific schools, project NSh-8081.2016.9.

\end{document}